# A consistent multi-user, multi-goal framework for assessing system performance with application to a sonar system

Dr. Colin M. Reed
formerly: QinetiQ Malvern
19 Broadlands Drive.
Malvern
Worcs
WR14 1PW
cm-p-reed@wr14-1pw.wanadoo.co.uk

Dr. Alan J. Fenwick
formerly: QinetiQ, Winfrith
now with: School of Engineering
University of Aberdeen
Aberdeen
AB24 3UE
ajfenwick@abdn.ac.uk

## Abstract

Agreeing suitability for purpose and procurement decisions depend on assessment of real or simulated performances of sonar systems against user requirements for particular scenarios. There may be multiple pertinent aspects of performance (e.g. detection, track estimation, identification/classification and cost) and multiple users (e.g. within picture compilation, threat assessment, resource allocation and intercept control tasks), each with different requirements. Further, the estimates of performances and the user requirements are likely to be uncertain. In such circumstances, how can we reliably assess and compare the effectiveness of candidate systems? This paper presents a general yet simple mathematical framework that achieves all of this. First, the general requirements of a satisfactory framework are outlined. Then, starting from a definition of a measure of effectiveness (MOE) based on set theory, the formulae for assessing performance in various applications are obtained. These include combined MOEs, multiple and possibly conflicting user requirements, multiple sources and types of performance data and different descriptions of uncertainty. Issues raised by implementation of the scheme within a simulator are discussed. Finally, it is shown how this approach to performance assessment is used to treat some challenging examples from sonar system assessment.

## 1. Introduction

It is becoming increasingly important to properly assess the performances of projected or in-service military surveillance systems either through simulation or trials. Efforts have been made to standardise a wide variety of pertinent scenarios (e.g. by the Studies Assumptions Group (SAG)) so that various sensing, tracking and identification systems, for example, can be compared on a 'level playing field' – at least by simulation. There are also efforts to standardise a suitable range of measures of performance (MOPs), which appear likely to extend into all military domains soon. Single Integrated Air Picture (SIAP) metrics are a good example as far as the air domain is concerned.

This paper is concerned with the next stage beyond MOPs, which we term measures of effectiveness (MOEs), a usage consistent with the original meaning of the term as discussed in the next section. The meaning of MOE has changed somewhat over the last decade. We take it to be the extent to which a MOP satisfies a declared user requirement (i.e. its original meaning). Involvement of the user requirement makes MOEs in this sense even more important than MOPs. Despite this, proposed MOEs and MOE systems have usually lacked a rigorous mathematical basis and that has raised concerns, not least from some military users. This paper aims to redress the perceived shortcomings of extant MOEs by presenting and justifying a rigorous mathematical framework for performance assessment that satisfies what are considered to be all the essential requirements. The framework aims to be generic, but is illustrated here by application to the assessment of a simulated sonar tracking system.

This paper is in three main parts. Section 2 discusses the background to system performance assessment. The third section proposes a framework for measures and justifies the claims made for it. The fourth and fifth sections collectively describe a multi-user multi-goal sonar tracking and identification system simulation. On a first reading, it is possible to skip the third section to see how the system performance is assessed within the framework.

## 2. Background

There are two prominent types of quantity for evaluating the performance of any system:

- measures of performance (MOPs);
- measures of effectiveness (MOEs).

A MOP is literally any convenient measure of what a system achieves in operation. Examples are measurement and estimation errors, report time delays, detection probability, mission cost, reliability and availability. A more comprehensive list for military command, control and communication (C3) systems is given by Sweet [6]; one for sensor management systems is given by Rothman and Bier [7]; and one for multi-sensor fusion systems is given by Llinas [8]. MOPs can be dimensional, can be expressed in different units and can have any range of values. The same system could have multiple MOPs to cover various aspects of performance.

In the sense used here, MOE expresses the extent to which a MOP satisfies a declared user requirement. As with MOPs, there could be multiple MOEs and there could be even more to cater for different users as well. A MOE is calculated from a MOP and a related user requirement. This meaning of MOE is taken from early studies of performance evaluation in the USA. Principal contributors include Levis at MIT [2], Matthis at Mitre Corporation [9] and Buede at George Mason University [10]. Further discussion of MOPs and MOEs can be found in the books by Waltz and Llinas [12] and Hall [13]. Despite the excellent foundation laid by these contributors, the term MOE has more recently become used, both in the USA and the UK, to indicate a MOP at war-fighting level, without reference to a user requirement [13], [23]. Our MOE, given the symbol $M$, is concerned with customer satisfaction or compliance. Also, we allow it to be used anywhere in a system and at any level, insofar as a customer has corresponding requirements.

Most formulae for calculating MOEs are heuristic and *ad hoc* and have led to concern in some quarters about the rationality, coherence and reliability of performance assessment. It is not our intention to 'name and shame' examples, some of which have served their users well, but we do list some deficiencies from a wider perspective as follows.

- MOEs with different physical units which cannot be compared or combined.
- MOE values in the range $[0, \infty]$ where it is difficult to appreciate the significance of values and their differences.
- Lack of a physical meaning.

- Limitation to one data type (usually continuous numeric).
- No provision for uncertainty in the measurements.
- Overriding user requirements by the laws of physics. That is a system is considered completely effective when it has achieved theoretical best performance, even though that performance falls well short of what the user has asked for.
- Overriding user requirements by the limitations of other systems. An example is not allowing the user to set measurement accuracy requirements beyond the precision allowed by the message protocol of a communication system.
- Threshold user requirements. For example, a track position error of 1.99 km is considered completely acceptable, while a slightly higher value of 2.01 km is considered completely unacceptable.

A more subtle example comes from a MOE formula reported by Blackman [14] in which the processing system is allowed to influence its own assessment. Consider a multi-target Kalman filter tracker using data from one sensor. Let $\varepsilon_j$ be the errors (between the true and predicted values for target azimuth, elevation and range) and let $C_{jj}$ be the Kalman filter variances. The reported MOE is given by:

$$M = \frac{1}{N}\sum_{i=1}^{N} S_i \quad ; \qquad \text{where} \quad S_i = \begin{cases} 1.0 \text{ for a confirmed track on target } i \text{ and } |d| \le d_t \\ 0.5 \text{ for a confirmed track and } |d| > d_t \\ 0.5 \text{ for a tentative track} \\ 0.0 \text{ for no track} \end{cases}$$

$$d^2 = \sum_{j=1}^{3} \frac{\varepsilon_j^2}{C_{jj}}$$

where $d_t$ is a threshold, and $d$ is a normalised error.

Several features of this MOE are reasonable. However, the normalised error involves Kalman filter variances, which are internal to the tracker. Thus, the tracker is being allowed to influence its own performance assessment. This is not reasonable when the quantities of operational importance are the errors $\varepsilon_j$.

The framework presented in this paper was developed by one of the authors starting around 1992. It was triggered by dissatisfaction of colleagues, including a seconded RAF officer, with existing MOE schemes for evaluation of air picture compilers. The intention was to develop a generic approach for MOEs from as rigorous a mathematical foundation as could be managed.

First, what were held to be the **essential** requirements for MOEs were elicited from colleagues working at Malvern on military C³I systems. These were as follows.

1. The MOEs should have a rigorous logical basis, i.e. not arbitrary in any way.
2. MOE values should have a suitable physical meaning and lie on a well-defined scale.
3. They should be easily calculated.
4. The framework should be able to handle any type of data (e.g. continuous numeric, discrete numeric, enumeration and multi-dimensional).
5. There should be a mechanism for combining MOEs to obtain an overall MOE.

6. The framework should be able to accept any expression of user requirements for system performance, whether objective or subjective in origin.

A MOEs framework that satisfied all of these requirements was successfully developed and used to assess the tracking performance of the UKADGE (UK Air Defence Ground Environment) system. The earliest generally available report containing a description of the framework was produced much later in 1996 [22]. A more up to date description, relevant to sonar system assessment, is provided in the next section.

## 3. The MOEs framework

### 3.1 Introduction

This section presents the key mathematical formulae for MOEs. These are the simplest possible, consistent with a coherent generic framework satisfying the six essential requirements listed at the end of section 2.

### 3.2 Set theory based definition of MOE

First, we consider here one of the simplest conceivable problems requiring a mathematical definition of MOE. It leads to a formula, which is the basis for all that follows.

The problem is stated in general, abstract (i.e. continuous set) terms. With reference to the Venn diagram of Figure 1, we start with a universal event space. Depending on the system output being evaluated, this represents the space of measurements or estimates (processed measurements). Henceforth, measurements and estimates will be embraced by the term 'observations'. It is supposed that a region of the space, with volume $V_o$, is uniformly occupied by observations from the system to be evaluated. Another region, with volume $V_s$, contains all observations that would be acceptable to a particular user, for example because they lay within certain error bounds. The problem is to derive the simplest mathematical definition of a MOE consistent with the six requirements in section 2.

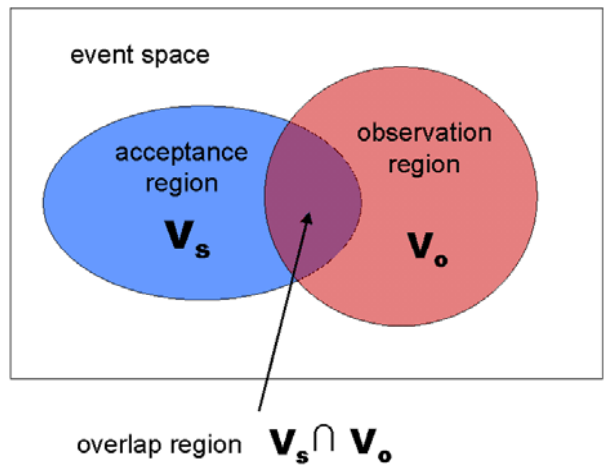

**Figure 1:** Venn diagram leading to simple set-theoretic definition of MOE.

It is shown in appendix A.1 that the result is:

$$M = \frac{V_o \cap V_s}{V_o} \qquad (1)$$

where $M$ denotes MOE. This notation is a shorthand for the ratio of the volume of the intersection (overlap) divided by the observation volume, and is dimensionless. It means 'the fraction of observations that are acceptable'. Values lie in the range [0,1]. Zero indicates complete non-compliance and unity indicates complete acceptance. These features satisfy requirements 1, 2 and 3.

**Remarks**

Another MOE candidate that deserves special mention is: $M = V_o \cap V_s / V_s$. This is capable of satisfying the requirements and it has the same complexity as result (1). However, it means 'the fraction of the user requirement that is satisfied by the observations', which is not the same as 'the fraction of observations that are acceptable'. As an example, this MOE would be appropriate for assessing coverage of a specified surveillance area by a sonar system. In this paper, further discussion is limited to the type of MOE expressed by result (1).

Set theory provides the most primitive basis for several important branches of mathematics (e.g. logic, arithmetic and probability theory). If a mathematical definition of MOE cannot handle a simple set-theoretic problem, like the one just treated, this is a reasonable ground for disqualifying it.

One can still obtain result (1) starting from discrete sets. The volume of a set is then taken to be proportional to the number of its unique elements and the intersection of volumes should be understood as the volume of the intersection of sets. Consequently, requirement 4 in section 2 is satisfied as well. We have not yet satisfied requirements 5 and 6, or shown that all the requirements can be satisfied for more general problems. All of this is done shortly.

### 3.3 Extension to uncertain data

To encompass more general problems, we must introduce uncertainty in both the observations and the user acceptance process. For present purposes, it is sufficient to deal with probabilistic uncertainty, but the formulation is capable of handling other models of uncertainty. With reference to Figure 2, we now have to treat probability density functions (PDFs), $\rho_o(\underline{x})$ and $\rho_s(\underline{x})$, for observations and user acceptance respectively. Quantity $\underline{x}$ is vector position in the space containing all possible observations and user specifications. In appendix A.2, it is shown that the generalisation of result (1) is the normalised overlap integral over all this space:

$$M = \int \frac{\rho_s(\underline{x})}{\rho_{s\max}} \rho_o(\underline{x}) \, d\underline{x} = \int f_s(\underline{x}) \, \rho_o(\underline{x}) \, d\underline{x} \qquad (2)$$

where $\rho_{s\max} = \max(\rho_s(\underline{x}))$. Quantity $f_s(\underline{x}) = \rho_s(\underline{x}) / \rho_{s\max}$ is called the **user function**. It is a measure of acceptance of an observation at $\underline{x}$, and takes values in the range [0,1]. As with MOE, zero indicates complete non-compliance and unity indicates complete acceptance. Indeed, if there is a single discrete observation with PDF $\delta(\underline{x} - \underline{x}')$, where $\delta$ denotes the Dirac delta function, result (2) becomes:

$$M = \int f_s(\underline{x})\,\delta(\underline{x} - \underline{x}')\,d\underline{x} = f_s(\underline{x}')$$  (3)

That is, the MOE equals the user function value. The user function encapsulates performance scoring in a general, yet simple manner. Therefore, it may be considered to satisfy requirement 6 in section 2. Practically, it can be expressed as a continuous parametric function (see the example given shortly) or as a discrete set of pairs, $\{\underline{x}_i, f_{si}\}$ for all $i$, (i.e. a look-up table).

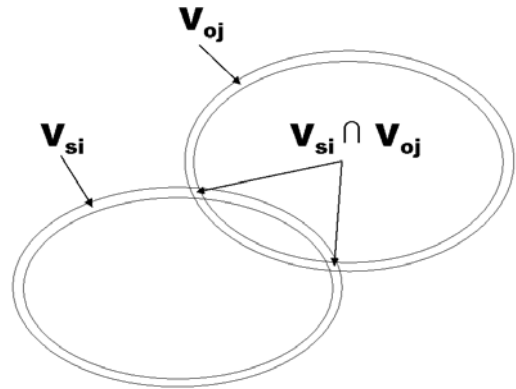

**Figure 2:** Intersecting iso-probability shells for acceptance and observation probability distributions.

Result (2) is an appropriately weighted version of 'the fraction of observations that are acceptable'. It contains result (1) as a special case, since the latter effectively has the uniform PDF:

$$\rho_k(\underline{x}) = \begin{cases} V_k^{-1} & \text{for all } k, \text{ if } \underline{x} \text{ is in } V_k \\ 0 & \text{otherwise} \end{cases}$$  (4)

**Example**

Suppose a system output has an unbiased error, $\varepsilon$, described by a 1-D Gaussian (Normal) PDF with standard deviation $\sigma_o$. Then, the observation PDF is:

$$\rho_o(\varepsilon) = \frac{1}{\sqrt{2\pi}\sigma_o} \exp\left(-\frac{\varepsilon^2}{2\sigma_o^2}\right)$$  (5)

The acceptability PDF, $\rho_s(\varepsilon)$, is also Gaussian, with standard deviation, $\sigma_s$. What is the MOE of the system over a large number of observations?

First, the user function is given by:

$$f_s(\varepsilon) = \frac{\rho_s(\varepsilon)}{\max(\rho_s(\varepsilon))} = \exp\left(-\frac{\varepsilon^2}{2\sigma_s^2}\right) \tag{6}$$

Finally, by result (2), the MOE is given by:

$$M = \int_{-\infty}^{\infty} f_s(\varepsilon)\,\rho_o(\varepsilon)\,d\varepsilon = \frac{1}{\sqrt{2\pi}\,\sigma_o} \int_{-\infty}^{\infty} \exp\left(-\frac{\varepsilon^2}{2}\left(\frac{1}{\sigma_o^2} + \frac{1}{\sigma_s^2}\right)\right) d\varepsilon \tag{7}$$

$$\therefore \quad M = \frac{\sigma_s}{\sqrt{\sigma_s^2 + \sigma_o^2}} \tag{8}$$

### 3.4 Relation to other work

Having derived some key results (i.e. (1) and (2)), this is the appropriate point at which to consider their relation to other work in this field.

Although obtained independently, result (1) is not original. As far as we are aware, it was first used by Bouthonnier and Levis [1] in 1984 as a basis for evaluation of Command and Control systems. It appears from a subsequent account by Levis [2] of the methodology they used that result (1) was simply taken as a definition of MOE, without derivation or justification. We have supplied the latter at least.

A MOE definition superficially similar to (2) has been used by Karam [3], but instead of our user function $f_s(\underline{x})$, he employs the **utility function** of von Neumann and Morgenstern [4]. Thus, Karam's definition is an expected utility. A more recent discussion of Karam's approach is given by Buede [5]. In general, there are important differences between user functions and utility functions.

- First, user function values are restricted to the range [0,1]. Utility functions (and expected utilities) need not be so limited – very large and negative values being possible. However, Karam's MOEs, like ours, are required to take values in the range [0,1]. This implies use of a special case of utility functions, which is indistinguishable from user functions as far as numerical range is concerned.

- Secondly, in contrast to utility functions, user functions are re-scaled PDFs, having values in the range [0,1]. Thus, our MOEs framework is an all-probabilistic one, whereas Karam's approach mixes two different quantities, probabilities and utilities. This is a crucial and advantageous difference. For example, in our approach, probability theory is already prescribed for problems such as combining different user functions. In contrast, it is often not obvious how to combine different utility functions; often a weighted sum is assumed and one is then faced with the problem of determining suitable weights. We would argue that the user function approach is more principled and easier to use in practice. This will become apparent in subsequent sections.

This paper is concerned with probabilistic uncertainty, but the framework can accommodate other models of uncertainty. In mostly unpublished work, this has been done for Dempster-Shafer, fuzzy and fuzzy rough set approaches. (See reference [22] for the fuzzy rough set definition of MOE.)

### 3.5 Special cases

Having established the general MOE formula (2) for probabilistic systems, we now show how it applies to some special cases.

**Case 1: Different observation and acceptability spaces**

Definition (2) applies where there is no information that is irrelevant to the user's purpose. However, when the process of observation automatically results in unusable information, the definition of probabilistic MOE can be extended. The form remains the same because the unwanted information can be integrated out to produce marginal distributions.

**Example**

Suppose we have the PDFs: $\rho_o(\underline{x}_1, \underline{x}_2)$ and $\rho_s(\underline{x}_1, \underline{x}_3)$. Then, quantity $Q$ is given by:

$$Q = \iiint \rho_s(\underline{x}_1, \underline{x}_3)\, \rho_o(\underline{x}_1, \underline{x}_2)\, d\underline{x}_1\, d\underline{x}_2\, d\underline{x}_3$$

$$\tag{9}$$

$$= \int \rho_s(\underline{x}_1)\, \rho_o(\underline{x}_1)\, d\underline{x}_1$$

since, by marginalisation: $\int \rho_o(\underline{x}_1, \underline{x}_2)\, d\underline{x}_2 = \rho_o(\underline{x}_1)$ and $\int \rho_s(\underline{x}_1, \underline{x}_3)\, d\underline{x}_3 = \rho_s(\underline{x}_1)$.
Quantity $Q_{max}$ is given by:

$$Q_{max} = \iiint \rho_s(\underline{x}_1, \underline{x}_3) \overbrace{\frac{\rho_o(\underline{x}_1, \underline{x}_2)}{\rho_o(\underline{x}_2)}}^{} \delta(\underline{x}_1 - \hat{\underline{x}}_1)\, d\underline{x}_1\, d\underline{x}_2\, d\underline{x}_3$$

$$= \iint \rho_s(\underline{x}_1, \underline{x}_3)\, \delta(\underline{x}_1 - \hat{\underline{x}}_1)\, d\underline{x}_1\, d\underline{x}_3 \quad \text{since} \quad \int \rho_o(\underline{x}_2)\, d\underline{x}_2 = 1$$

$$= \int \rho_s(\underline{x}_1)\, \delta(\underline{x}_1 - \hat{\underline{x}}_1)\, d\underline{x}_1 = \rho_s(\hat{\underline{x}}_1) = \max(\rho_s(\underline{x}_1)) \tag{10}$$

Note that only the observation probability mass over $\underline{x}_1$ is concentrated (at $\hat{\underline{x}}_1$), because it is the only variable in the observation PDF that also occurs in the acceptability PDF. Finally:

$$M = \frac{Q}{Q_{max}} = \int \frac{\rho_s(\underline{x}_1)}{\max(\rho_s(\underline{x}_1))}\, \rho_o(\underline{x}_1)\, d\underline{x}_1 = \int f_s(\underline{x}_1)\, \rho_o(\underline{x}_1)\, d\underline{x}_1 \tag{11}$$

This is reasonable behaviour for MOEs. The probability mass of observations over irrelevant variables is aggregated onto the relevant ones. Similarly, the probability mass of acceptability over hidden (unobserved) variables is aggregated onto the observed ones. However, one should always be careful to note which variables are encompassed by the resulting MOE after marginalisation.

**Case 2: Discrete data types**

Result (2) is stated for continuous variables, but it handles all the data types listed in requirement 4 of section 2 as special cases. First, consider discrete numeric data. Here, observations occur as a set of pairs, $\{\underline{x}_i, P_{oi}\}$ for $i$ = 1 to $I$ say, where $P_{oi} = P_o(\underline{x}_i)$ is the probability of observation $\underline{x}_i$. The observation PDF is:

$$\rho_o(\underline{x}) = \sum_{i=1}^{I} P_{oi} \; \delta(\underline{x} - \underline{x}_i) \tag{12}$$

and the corresponding acceptability PDF is:

$$\rho_s(\underline{x}) = \sum_{i=1}^{I} P_{si} \; \delta(\underline{x} - \underline{x}_i) \tag{13}$$

where the acceptance probabilities $\{P_{si}\}$ are set by the user. Inserting these two PDFs into result (2) leads to the MOE for discrete numeric quantities:

$$M = \frac{1}{P_{s\,\max}} \sum_{i=1}^{I} P_{si} \; P_{oi} = \frac{\mathbf{P}_s^T \; \mathbf{P}_o}{\max(\mathbf{P}_s)} = \mathbf{f}_s^T \; \mathbf{P}_o \tag{14}$$

This is simply the scalar product of user function and observation probability vectors, $\mathbf{f}_s = \begin{pmatrix} f_{s1} & f_{s2} & \cdots & f_{sI} \end{pmatrix}^T$ and $\mathbf{P}_o = \begin{pmatrix} P_{01} & P_{o2} & \cdots & P_{oI} \end{pmatrix}^T$ respectively, where superscript $T$ denotes matrix transposition. It is also a projection of the user function vector onto the observation probability vector.

Since index $i$, denoting vector component, serves only as a label, it can equally well be non-numeric. Therefore, result (14) handles enumeration data, where the observations occur as the set of pairs $\{i, P_{oi}\}$, $i$ being a class label for example.

**Example**

Suppose a naval user requires that a surveillance system correctly classify an enemy vessel. The system outputs an identity vector estimate, in the form of probabilities over three possible classes: $E$ - enemy, $N$ - neutral, and $F$ - friend. This has the value: $\mathbf{P}_o = \begin{pmatrix} P_E & P_N & P_F \end{pmatrix}^T = \begin{pmatrix} 0.60 & 0.25 & 0.15 \end{pmatrix}^T$. The user function vector for this situation is: $\mathbf{f_s} = \begin{pmatrix} f_{sE} & f_{sN} & f_{sF} \end{pmatrix}^T = \begin{pmatrix} 1.0 & 0.0 & 0.0 \end{pmatrix}^T$, since only that part of the estimate in favour of the true enemy class is acceptable. By result (14), the MOE for the estimate is:

$$M = \mathbf{f}_s^T \; \mathbf{P}_o = \begin{pmatrix} 1.0 & 0.0 & 0.0 \end{pmatrix} \begin{pmatrix} 0.60 \\ 0.25 \\ 0.15 \end{pmatrix} = 0.60 \tag{15}$$

**Case 3: Uniform user functions**

If the user function is uniform over a finite region, result (2) reduces to the bounded volume integral:

$$M = \int_{V_s} \rho_o(\underline{x}) \, d\underline{x} \tag{16}$$

This can be recognised as a multi-dimensional version of a familiar acceptance testing formula. To illustrate this, suppose we have a measuring device whose error $\varepsilon$ follows a 1-D Gaussian PDF with standard deviation $\sigma_o$. A potential user of the device requires errors to be within the bounds $\pm \Delta \varepsilon$. By result (16), the MOE is:

$$M = \frac{1}{\sqrt{2\pi}\sigma_o} \int_{-\Delta\varepsilon}^{\Delta\varepsilon} \exp\left(-\frac{\varepsilon^2}{2\sigma_o^2}\right) d\varepsilon = \mathrm{erf}\left(\frac{\Delta\varepsilon}{\sqrt{2}\sigma_o}\right) \tag{17}$$

where 'erf' denotes the Error Function. As usual, the MOE is the fraction of measurements which are acceptable and, depending on how stringent the user's requirements are, (s)he may accept $M$ = 0.95 or 0.99 typically, these corresponding to 'confidence levels'.

In the majority of situations, it is inappropriate to have a user function with a hard cut-off, as defined by the boundary of volume $V_s$ here.

**Case 4: Sampled observations**

In practice, most observations are samples, so this special case leads to the most useful of all MOE formulae. Samples take the form of a set of discrete values, $\{\underline{x}_i\}$ for $i$ = 1 to $I$ say. The associated probabilities do not appear because all the values implicitly have the same probability of $1/I$. This is consistent with either random sampling **or** a state of ignorance. The corresponding PDF is:

$$\rho_o(\underline{x}) = \frac{1}{I} \sum_{i=1}^{I} \delta(\underline{x} - \underline{x}_i) \tag{18}$$

Insertion of this into result (2) produces:

$$M = \frac{1}{I} \sum_{i=1}^{I} f_s(\underline{x}_i) = \overline{f_s} \tag{19}$$

The MOE is simply the mean of the user functions evaluated over the samples. Many measurements come as just single samples; the MOEs are then just the corresponding user function values.

**3.6   Multiple observers and users**

The set-theoretic problem examined in section 3.2, which led to a mathematical definition of MOE, prompts the question of what happens when there are multiple observation regions (corresponding to similar outputs from different systems) and multiple acceptability regions (corresponding to different users, who may not agree on their requirements)? It is emphasised that all of these share the same observation space. In these circumstances, result (2) simply generalises to:

$$M = \int F_s\left(f_{s1}, f_{s2}, \cdots, f_{sI}\right) F_o\left(\rho_{o1}, \rho_{o2}, \cdots, \rho_{oJ}\right) d\underline{x} \qquad (20)$$

where for brevity the arguments of the $I$ user functions and $J$ observation PDFs are omitted. In the absence of further information, all observations are treated identically, as are all users. Then, the separate combination functions, $F_s$ and $F_o$, while generally different, are always symmetric (i.e. showing no preference for the order of their arguments). $F_s$ produces a simple user function and $F_o$ produces a simple PDF. The combination functions are elaborated as follows:

- **Combining observations.** The basic requirements of symmetry and producing a simple PDF are common to those for suitable data fusion approaches and these are what would be applied here. Even within probabilistic systems, there is more than one way of combining PDFs, for example depending on the availability of priors and whether the PDFs are independent or not. If we are given just the separate PDFs and no further information, there is no option but to use the combination law for independent PDFs:

$$F_o\left(\underline{x}\right) = \frac{\prod_{j=1}^{J} \rho_{oj}\left(\underline{x}\right)}{\int \prod_{j=1}^{J} \rho_{oj}\left(\underline{x}\right) d\underline{x}} \qquad (21)$$

If there is dependence between the observation sources, result (21) will not be correct, but it may serve as a reasonable approximation. Any MOE based on it will have to be qualified accordingly. Although result (20) can handle combination of observations, this is better achieved with more complete information before the performance evaluation stage.

- **Combining user functions.** User functions are re-scaled PDFs of the form: $f_{si}\left(\underline{x}\right) = \rho_{si}\left(\underline{x}\right)/\max\left(\rho_{si}\left(\underline{x}\right)\right) \equiv \rho\left(\underline{x} \mid \underline{S}_i\right)/\max\left(\rho\left(\underline{x} \mid \underline{S}_i\right)\right)$, where $\underline{S}_i$ is a set of parameter values determining the form of the acceptability PDF for user '$i$'. Some of those values will be extracted from a continuum of possible values. It turns out that we can combine these functions in a potentially wide variety of ways, reflecting how we might wish to satisfy the users as a whole. For $I$ users, these are generated by the symmetric transform:

$$\sum_{k=1}^{I} \binom{I}{k} F_{sk}^{k} \, Z^{I-k} = \prod_{i=1}^{I}\left(Z + f_{si}\right) - Z^{I} \qquad (22)$$

The derivation of this result is outlined in appendix A.3. The various combination functions, $F_{sk}$, are found by comparing similar powers of the dummy variable $Z$. For the example of four users ($I = 4$), this produces:

$$F_{s1} = \frac{1}{4}\left(f_{s1} + f_{s2} + f_{s3} + f_{s4}\right)$$

$$F_{s2} = \left(\frac{1}{6}\left(f_{s1}\, f_{s2} + f_{s1}\, f_{s3} + f_{s1}\, f_{s4} + f_{s2}\, f_{s3} + f_{s2}\, f_{s4} + f_{s3}\, f_{s4}\right)\right)^{\frac{1}{2}}$$

$$F_{s3} = \left( \frac{1}{4} \left( f_{s1} \, f_{s2} \, f_{s3} + f_{s1} \, f_{s2} \, f_{s4} + f_{s1} \, f_{s3} \, f_{s4} + f_{s2} \, f_{s3} \, f_{s4} \right) \right)^{\frac{1}{3}}$$

$$F_{s4} = \left( f_{s1} \, f_{s2} \, f_{s3} \, f_{s4} \right)^{\frac{1}{4}} \tag{23}$$

Generally, there are as many combination functions as users. Each $F_{sk}$ has the form of a simple user function with values in the range [0,1]. When all the user functions are identical ($f_{si} = f_s$ for all $i$), so are all the combination functions ($F_{sk} = f_s$ for all $k$). $F_{s1}$ is the most accommodating combination function, since satisfying any of the users makes a contribution to overall acceptability. $F_{s2}$ is more demanding, requiring that pairs of users be jointly satisfied and so on. The most demanding function, $F_{s4}$, requires that all the users be jointly satisfied.

So far, we have treated all users or their equivalent n-tuplets without distinction. However, we can allow legitimately for relative importance by weighting. For example, $F_{s1}$, would take the form:

$$F_{s1} = W_1 \, f_{s1} + W_2 \, f_{s2} + W_3 \, f_{s3} + W_4 \, f_{s4} \tag{24}$$

where the weights, $W_i$, are positive and sum to unity. Effectively the user functions are re-scaled. In the case of $F_{s4}$ above, the same formula stands irrespective of importance weighting.

### 3.7 Combining MOEs

So far, we have shown that requirements 1 to 4 and 6 in section 2 can be satisfied when dealing with general probabilistic systems. We now address the remaining requirement 5, concerned with combining different MOEs to produce an overall MOE.

First, it should be noted that all MOE values in this framework mean exactly the same thing irrespective of what they physically relate to. A MOE means the extent to which a user is satisfied by an observation. Therefore, it is reasonable to combine MOEs directly to obtain an overall MOE. In practice, it is unlikely that a user with multiple requirements will do other than specify separate (i.e. independent) user functions. Observations over the multiple variables involved may show dependence or independence. Therefore, it is appropriate to have MOE combination formulae for two cases as follows:

- **Independent observations**

For this case, the observation PDF over multiple variables $x_1$ to $x_J$ has the separated form, $\rho_o(\underline{x}) = \rho_{o1}(x_1) \, \rho_{o2}(x_2) \cdots \rho_{oJ}(x_J)$, and the user function has a similar form, $f_s(\underline{x}) = f_{s1}(x_1) \, f_{s2}(x_2) \cdots f_{sJ}(x_J)$. Inserting these into the basic MOE formula (2) and carrying out the integration produces the overall MOE:

$$M = \prod_{j=1}^{J} M_j \quad \text{where} \quad M_j = \int f_{sj}(x_j)\, \rho_{oj}(x_j)\, dx_j \qquad (25)$$

Because of its multi-dimensional character it can be considered a MOE volume, still taking values in the range [0,1]. However, it becomes difficult to appreciate what those values mean. For example, if we have 10 variables and the MOE for each is 0.80, then, by result (25), the overall MOE is about 0.11. This is a good value, but instinctively we would compare it with a 1-D MOE and think it poor. To refer all MOEs to a common 1-D basis, we should take instead of result (25) the geometric mean:

$$M^{'} = \left( \prod_{j=1}^{J} M_j \right)^{\frac{1}{J}} \qquad (26)$$

This ensures that when all the MOEs are identical, the overall MOE is equal to one of them.

Result (26) is quite different from a widely used utility equivalent involving a weighted sum of separate utilities. However, each has a parallel with results in section 3.5. Result (26) is consistent with having to satisfy a user on all points. The utility result allows us to get away with satisfying some of them well. Arguably, result (26) is the better basis for performance assessment.

- **Dependent observations**

When the observations are dependent, the observation PDF is no longer a product of one-dimensional PDFs. It is possible to develop an algebraic formula for overall MOE from this, even including dependent user functions, but it is complicated and almost unusable in practice. It is suggested that an approach to assessment in this case is to take a large number $K$ of independent random samples $\{x_{1k}, x_{2k}, \cdots, x_{Jk}\}$ for $k =$ 1 to $K$ from the observation PDF. Then, by result (19) the (1-D referred) overall MOE is closely approximated by:

$$M^{'} = \left( \frac{1}{K} \sum_{k=1}^{K} \left( f_{s1}(x_{1k})\, f_{s2}(x_{2k}) \cdots f_{sJ}(x_{Jk}) \right) \right)^{\frac{1}{J}} \qquad (27)$$

It transpires that the geometric mean is appropriate for independent user functions. This results from the observation PDF taking a separated form when the observations are discrete. The situation with dependent user functions is not clear.

### 3.8 Allowance for 'uncertain ground truth'

In many applications the quantity we wish to assess is measurement error. When the system is simulated, the 'ground truth' is conveniently available, though sometimes it may be a challenge to associate it with the corresponding observation. The situation is more difficult with live trials. An example is to establish exactly where a target vessel is located for assessing the accuracy of a tracking system. Sometimes a more accurate, independent measurement of location can be arranged, for example through deployment of a co-operative target vessel fitted with GPS. If none of the errors is negligible how can we calculate the MOE of the tracking system? This can be done in a rigorous manner, even for single-sample measurements, provided the error distribution of the independent measurement is known.

If $\underline{x}'$ is a sampled measurement from the system being evaluated and $\underline{y}'$ is a corresponding measurement from the more accurate system, with error PDF $\rho_\lambda(\lambda)$, the MOE is given by:

$$M = \int f_s(\underline{\varepsilon})\, \rho_\lambda\left(\underline{\varepsilon} - \underline{x}' + \underline{y}'\right) d\underline{\varepsilon} \tag{28}$$

where $\underline{\varepsilon}$ is the measurement error for the system being evaluated. We never need to know this error, but we must have a user function for it. The derivation of this result is given in appendix A.4

**Example 1:** If the ground truth is known accurately, $\rho_\lambda(\lambda) = \delta(\lambda)$ and result (28) becomes:

$$M = \int f_s(\underline{\varepsilon})\, \delta\left(\underline{\varepsilon} - \underline{x}' + \underline{y}'\right) d\underline{\varepsilon} = f_s\left(\underline{x}' - \underline{y}'\right) = f_s\left(\underline{\varepsilon}'\right) \tag{29}$$

since $\underline{y}'$ now coincides with the ground truth. This is the usual, single-sample MOE result, albeit in terms of error $\underline{\varepsilon}'$.

**Example 2:** Here we consider a multi-dimensional example with a Gaussian error PDF as follows:

$$\rho_\lambda(\underline{\lambda}) = \frac{1}{\sqrt{\det(2\pi C_\lambda)}} \exp\left(-\frac{1}{2} \underline{\lambda}^T C_\lambda^{-1} \underline{\lambda}\right) \tag{30}$$

'$\det$' indicates 'determinant of' and $C_\lambda$ is the covariance matrix of the error vector $\underline{\lambda}$. A Gaussian user function is taken for error vector $\underline{\varepsilon}$ as follows:

$$f_s(\underline{\varepsilon}) = \exp\left(-\frac{1}{2} \underline{\varepsilon}^T C_s^{-1} \underline{\varepsilon}\right) \tag{31}$$

In practice, the covariance matrix $C_s$ is likely to be diagonal to avoid specifying non-zero off-diagonal elements. Applying result (28) produces:

$$M = \int \frac{1}{\sqrt{\det(2\pi C_\lambda)}} \exp\left(-\frac{1}{2}\left(\underline{\varepsilon}^T C_s \underline{\varepsilon} + \left(\underline{\varepsilon} - \underline{x}' + \underline{y}'\right)^T C_\lambda^{-1}\left(\underline{\varepsilon} - \underline{x}' + \underline{y}'\right)\right)\right) d\underline{\varepsilon} \tag{32}$$

After some matrix manipulation and integration, we finally obtain:

$$M = \sqrt{\frac{\det(C)}{\det(C_\lambda)}} \exp\left(-\frac{1}{2}\left(\underline{x}' - \underline{y}'\right)^T C_o^{-1}\left(\underline{x}' - \underline{y}'\right)\right) \tag{33}$$

where:

$$C^{-1} = C_s^{-1} + C_\lambda^{-1} \qquad \text{and} \qquad C_o^{-1} = C_\lambda^{-1} - C_\lambda^{-1} C C_\lambda^{-1} \qquad (34)$$

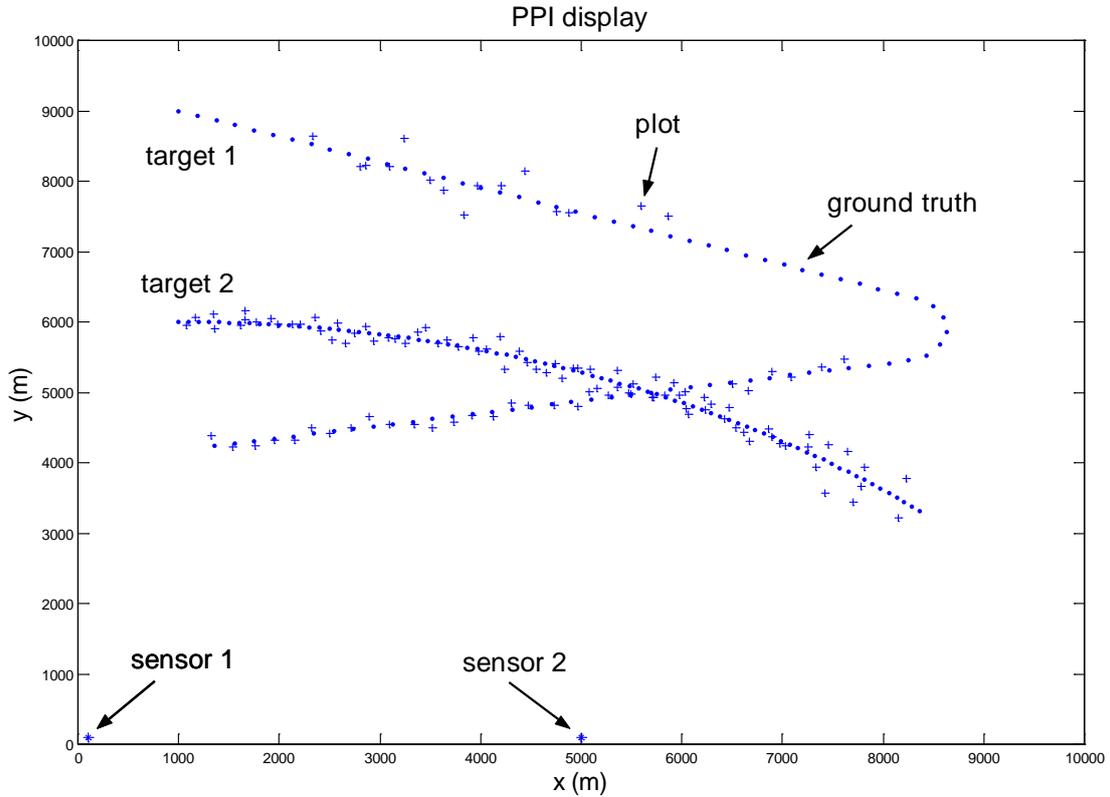

**Figure 3:** Plan position indicator display showing the (fixed) locations of two passive sonar sensors, the trajectories of two submarine targets (dots) and the plots from triangulation of sensor bearing measurements (+ signs). The targets start moving to the East from their labels. Target 1 is a friend and target 2 is an enemy.

## 4. Simulation

We have now described a framework for system performance assessment which is able to deal consistently with multiple goals, multiple users with different, even conflicting, requirements, and uncertainty The use of the measures/metrics is now demonstrated for a sonar tracking problem as follows.

- The scenario features two submarines of different allegiances (enemy and friend) following trajectories as shown in Figure 3.
- There are two stationary passive sonar sensors, each with uncertain detection, bearing estimates and identity estimates (in the form of identity vectors, as defined in STANAG 4162 [19]). The sonar sensors might be the 'best' two in a sonobuoy field, for example.
- Time-synchronised bearing measurements on each target are triangulated and plan-position estimates (plots) with associated covariance matrices calculated. Identity vectors are fused according to a conservative arithmetic averaging law.
- A simple Kalman-filter tracker, with perfect plot-to-track association through the use of identity estimates, converts the plots to tracks. There are two variants of this tracker, the difference being in just the levels of system (plant) noise employed in the state transition relation (one of the key tracking relations). The system noise determines the smoothness and responsiveness of tracking.

- Finally, the track and ground-truth data are accessed to calculate various MOEs as described in the next section. We cater for two (friendly) naval users with somewhat different requirements and, therefore, user functions.

Thus, the problem features a minimal multi-target scenario, different trackers whose multi-criteria performances are to be evaluated and compared, and users with different requirements.

In order to obtain unclassified sonar track information on the submarine targets together with ground truth information, all of which are needed to fully test the capabilities of the MOEs framework, we have had to develop a new simulator/test-bed. It is programmed in Matlab and runs on a PC. Ground-truth data produced by the scenario generator includes target times, positions, speeds, headings and identities. These are stored in an ASCII file that is accessed by the sonar sensor simulators. Low-fidelity modelling of fictitious sensors is employed; this is sufficient for demonstration purposes. The plot information from sensor-level fusion is stored in an ASCII file that can be accessed by the trackers. The resulting tracks are stored in ASCII files. It may be noted that there are better types of tracking filter than the simple Kalman one for this application, but the latter is sufficient for demonstration. A data flow diagram for the test-bed is shown in Figure 4

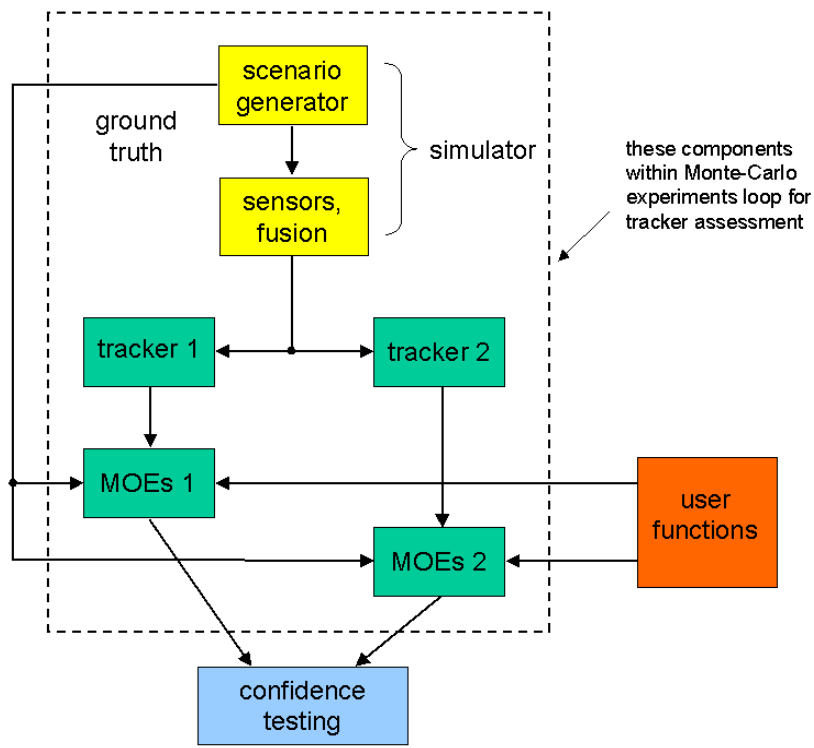

**Figure 4:** Schematic diagram of test-bed.

In our simulator, it was straightforward to associate the tracks with the corresponding ground truth. This was because the sonar reports were labelled with the targets of origin, the correct report was associated with each track (on the basis of identification information), and the labels were carried through. Also, there was time-synchronisation between ground-truth data, sonar reports and tracks. This is the ideal situation, if it can be arranged, because track-to-ground truth association is correct and certain.

In the case of trackers of the PDA (Probabilistic Data Association) type [17], (which were not tested here), plots from different targets may be used to update the same track. One cannot use plot-to-track association probabilities, say, to help calculate the proper MOE, because that would be using the tracker to assist in its own evaluation – rather like depending on self-assessment. A strategy avoiding this is to equally weight all target-to-track associations, i.e. the complete ignorance position, and calculate the corresponding MOE as follows:

$$M_k = \frac{1}{N} \sum_{i=1}^{N} M_{ki} \qquad (35)$$

where $M_{ki}$ is the MOE for track $k$ against ground truth/target $i$, there being $N$ such associations. This is favoured and has been used by one of the authors in other work [21]. However, it will be over-conservative in its MOE evaluations when targets are close to each other.

Another approach, promoted by Drummond [20] in the USA, is to independently associate ground truth with tracks using a global nearest neighbour, or better still a multi-frame assignment, method. The former finds the most likely pairings in a single time frame using, for example, a Munkres type algorithm [18]. The latter does the same thing over a small number of time frames and produces more reliable results. The approach is likely to produce accurate MOE results and avoids the need for passing target labels through the simulated sensors and trackers. The latter is an important consideration when modifying existing software is not a practical proposition. However, the assignment is not guaranteed to be correct and there may be no indication of when it has failed. Although failure may be infrequent, it can cast suspicion on all the MOE results. A completely satisfactory, general solution to the ground truth-to-track association problem is still awaited.

## 5. Results

This section presents some example results demonstrating most of the capability of the MOEs framework described earlier.

Some important background information is summarised first. The two submarine targets in the scenario have the (true) allegiances:

| Target identity no. | Allegiance |
|---|---|
| 1 | F – friend |
| 2 | E – enemy |

User functions are characterised as follows. User 1 chooses a Gaussian exponential function for errors in position, speed and heading (see equation (6)). He sets the value of the single parameter involved (standard deviation $\sigma_s$) for each of these variables. Generally, when $x = \sigma_s$, $f_s(x) = 0.6065$. User 2 prefers the user function:

$$f_s(x) = \frac{a^2}{(a^2 + x^2)} \qquad (36)$$

In this case, the single parameter $a$ is the value of $x$ at which $f_s(x) = 0.5$.

The user-specified parameter values are tabulated as follows:

| Variable | Units | User 1 - $\sigma_s$ | User 2 - $a$ |
|---|---|---|---|
| Position error | M | 0.5 | 0.2 |
| Speed error | m/s | 2.0 | 1.0 |
| Heading error | deg. | 5.0 | 2.0 |
| Identity | probability 3-vector | [0.2,0.2,0.6] for F | [0.0,0.0,1.0] for F |
| | | [0.6,0.2,0.2] for E | [[1.0,0.0,0.0] for E |

Examples of single MOE (for each variable) versus time plots, for tracks on target 1 from tracker 1 and for user 1, are shown in Figure 5. When the target is not detected, none of these MOEs can be calculated; hence the absence of values over some time instants (corresponding to sonar/tracker time frames). Even where they do exist, MOE values can be highly variable depending on what is happening in the scenario. The reader is reminded that these MOEs tell us directly the extent to which user requirements are being met at each time instant.

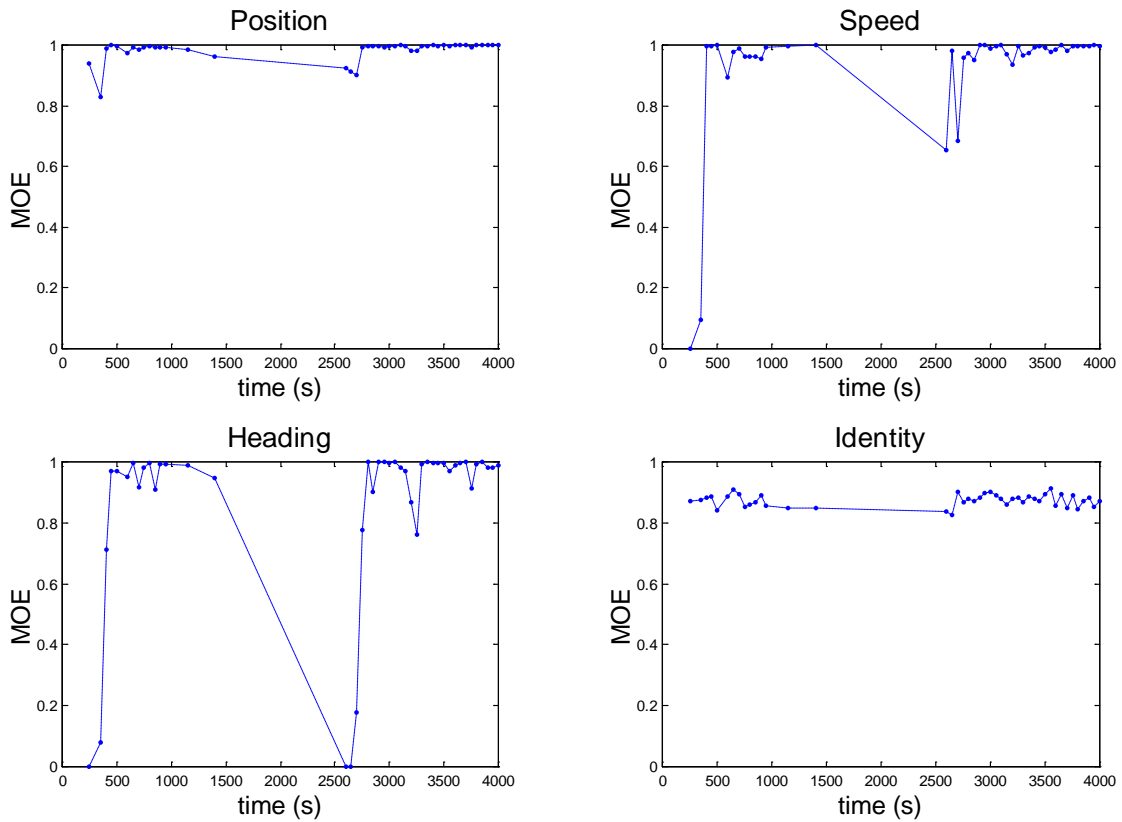

**Figure 5:** Separate MOEs versus time for tracker 1 and user 1. The MOEs are for errors in position, speed and heading, and for identity. The dotted line connecting points is not meant to interpolate MOE values (which may not exist between sample times).

Access to single MOE versus time plots allows the most detailed analysis and performance assessment. However, for some purposes, it may be preferable to work with overall MOEs. In the above scenario, there are at least three ways in which MOEs can be combined, producing:

- a multi-target MOE for each variable. Logically, this entails an arithmetic mean of the single target MOEs available at each time instant;
- a multi-variable MOE for each target. This entails a geometric mean of the single variable MOEs available, as discussed in section 3. Figures 6a and 6b show such combined MOEs for targets 1 and 2 respectively, both for tracker 1 and user 1;
- a multi-target, multi-variable MOE, entailing synthesis of the above, i.e. an arithmetic mean over targets of geometric means over variables for each target. (One has to be careful about what is intended, as there is also a multi-variable, multi-target MOE involving a geometric mean of an arithmetic mean, which is not the same.) Figure 6c shows the results of the synthesis intended, using the results in Figures 6a and 6b.

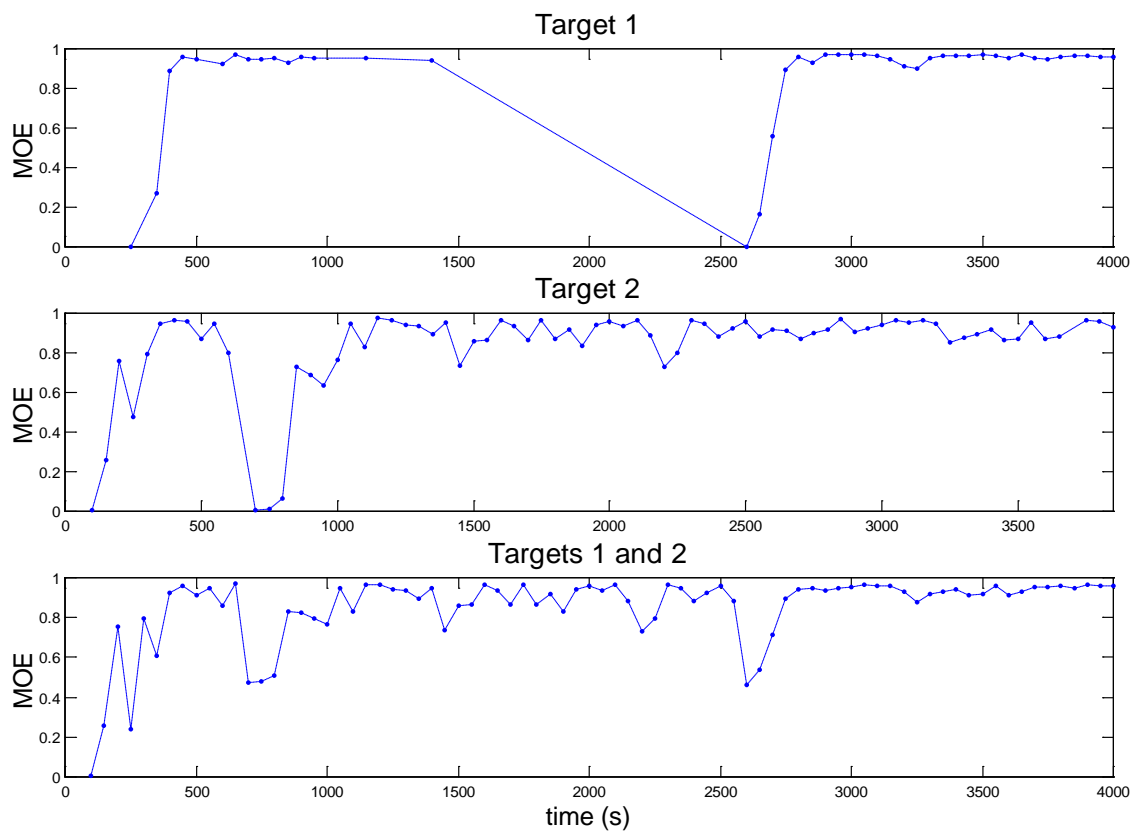

**Figure 6:** Overall MOEs for target 1 (6a), target 2 (6b), and both targets (6c). In each case the overall MOE covers position, speed, heading and identity variables. The results are for tracker 1 and user 1.

MOEs for combined users are shown in Figure 7. First, in Figures 7a and 7b, we have overall (multi-target, multi-variable) MOE versus time plots for users 1 and 2 respectively. (Figure 7a repeats Figure 6c for convenience.) Figure 7c shows the effect of combining the user functions according to an (equally weighted) arithmetic mean. This is the least demanding way in which the users could be satisfied collectively. In contrast, Figure 7d shows the effect of combining the user functions according to a geometric average. This is the most demanding way and that is reflected in the generally lower MOE values, though in this case the differences are not great.

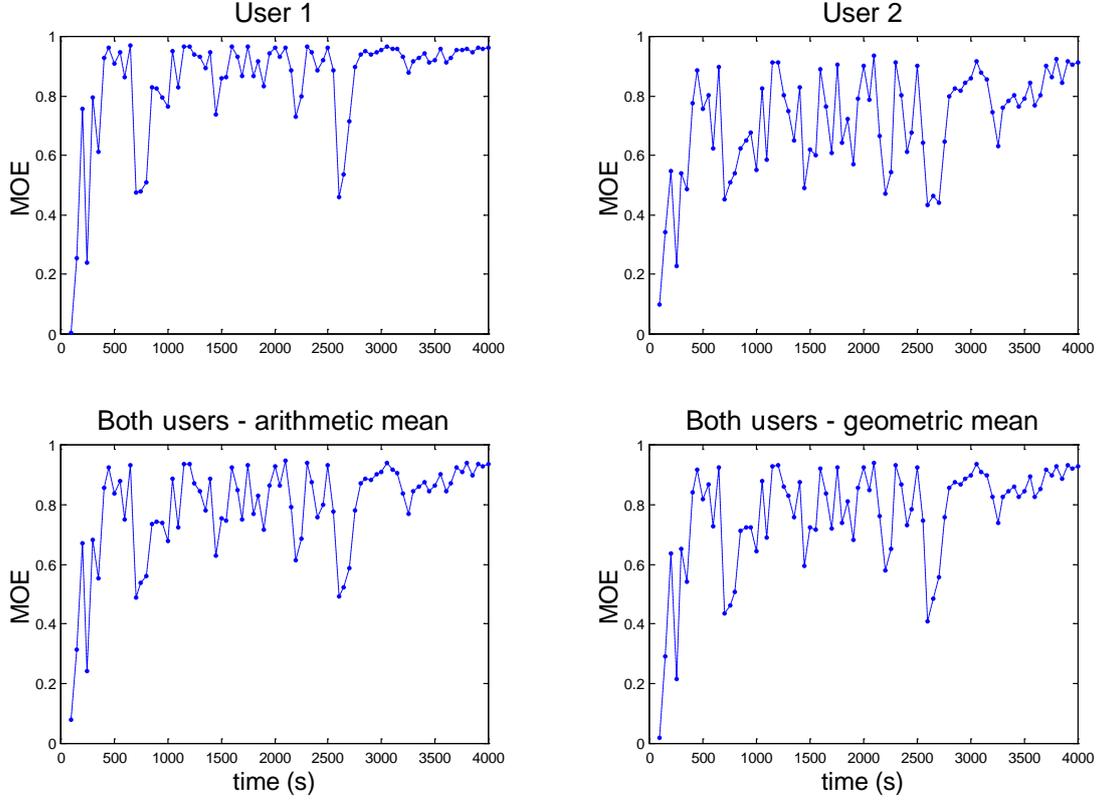

**Figure 7:** Overall MOEs (variables and targets) for tracker 1. The plots are for user 1 (7a), user 2 (7b), an (conservative) arithmetic mean of user functions (7c) and finally a (stringent) geometric mean of user functions (7d).

Finally, Figure 8 shows the results of a statistical comparison between overall (multi-target, multi-variable) MOEs from trackers 1 and 2. Here, a geometric mean of user functions is taken. The results from 20 Monte-Carlo simulation runs have been analysed to provide the difference in mean MOE values and the 95% confidence limits at each time instant. A difference of means exceeding the confidence limits indicates a significant difference between the performances of the trackers at this time instance. Following standard statistical theory [15], the confidence limits are given by:

$$\Delta M = t \sqrt{\left(\frac{1}{n_1} + \frac{1}{n_1}\right)\left(\frac{n_1 V_1 + n_2 V_2}{n_1 + n_2 - 2}\right)} \tag{37}$$

where $n_1$ and $n_2$ are the numbers of MOE samples for trackers 1 and 2 respectively, $V_1$ and $V_2$ are the corresponding sample variances, and $t$ is the Student **t** value for 95% confidence and $(n_1 + n_2 - 2)$ degrees of freedom. Strictly, result (37) assumes an underlying Gaussian distribution for MOE values. That cannot be the case here, since, at the very least, MOE values are restricted to the range [0,1]. However, if the numbers of samples are reasonably large, we can appeal to the Central Limit Theorem [16] to conclude that result (37) should be approximately true.

The results in this section illustrate a wide range of the MOE analyses that can be made to show a customer directly how a sonar tracking system meets his requirements and how it compares with other such systems. The results are scenario specific. There are many other MOEs not considered here, e.g. for probability of detection, track breaking, false tracks, purchase and running costs, averages over the scenario duration and averages over multiple scenarios. All of these can be handled easily within the same overall framework.

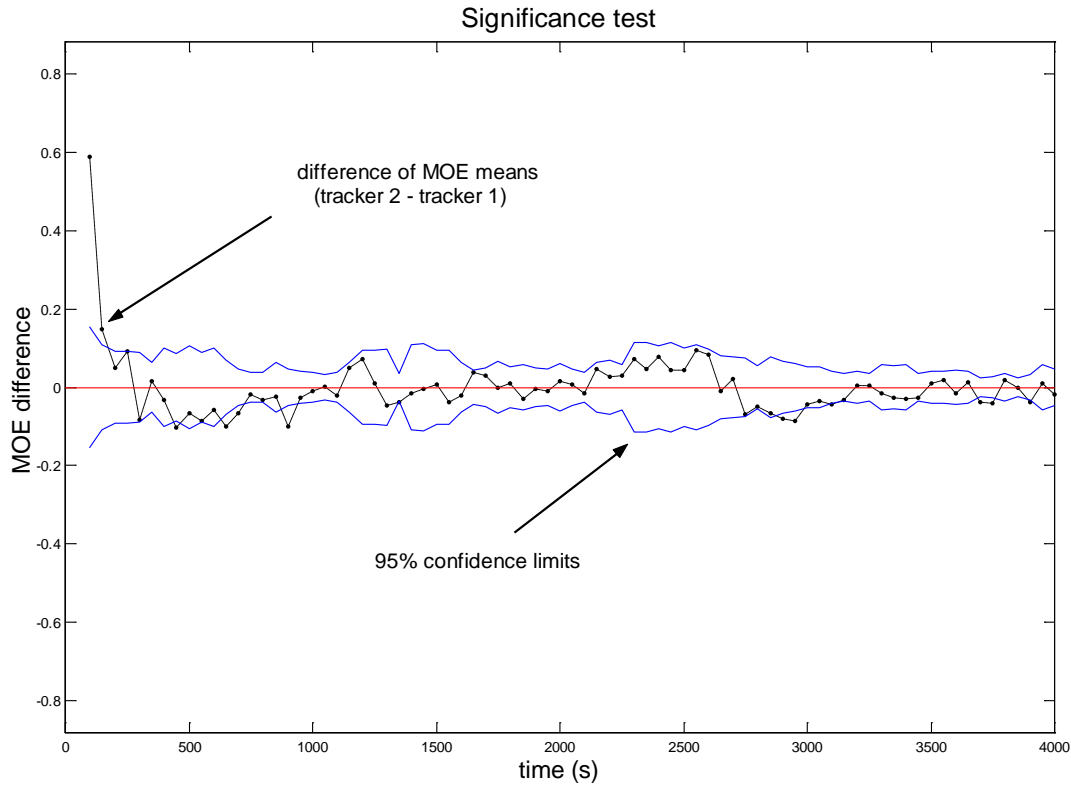

**Figure 8:** Significance testing of the difference in overall MOEs (over variables, targets and users) from trackers 1 and 2. Of the 79 points, 13 (i.e. about 16%) fall outside the confidence limits. This is about 3 times more than the 4 points (i.e. 5%) that can be expected from random variation. On the whole, tracker 2 is significantly less effective than tracker 1.

## 6. Conclusions

We have developed a mathematically rigorous, yet minimally simple, framework for calculating measures of effectiveness (MOEs). The basic mathematical definition of our MOE means *the fraction of observations that meet user requirements* (allowing for appropriate weightings). MOE values are easy to understand and calculate. The definition of our MOE has an all-probabilistic basis (in contrast to utility theory, for example) which enables it to naturally and correctly handle:

- multi-variable performance assessment of any system against user requirements;
- multiple object (e.g. target) scenarios;

- multiple users who must be satisfied according to different strategies;
- multiple sources of the performance data to be assessed.

Most of the above has been demonstrated in this paper for a simulated sonar tracking application. This has featured two submarine targets, two passive sonar sensors whose measurements are fused to provide plan-position plots, two simple Kalman filter trackers with different system (plant) noise settings whose performances are to be compared, and two users with different requirements. In all cases, some of them quite involved, it has proved quite straightforward to correctly calculate and combine user requirements and MOEs. Further, it has been shown how the MOEs for different trackers can be submitted to conventional confidence testing to establish the significance of differences.

This represents a fairly extensive testing and comparison of sonar tracker performances directly in terms of effectiveness against user requirements – the quantity of prime importance for procurement decisions. It demonstrates that the framework provides a powerful tool for general performance assessment of simple systems and is capable of dealing also with complex systems.

In practice, establishing user requirements, especially at lower system levels (e.g. sensors, trackers, classifiers, picture compilers and threat assessment) is not trivial – though this has been done for the UK air defence domain, for example. Methodologies for derivation of quantitative user requirements would be of great value. The MOE framework described here could assist with that. A principle of MOE balancing can be invoked, that theoretically allows derivation of lower-level system requirements from higher-level ones - given knowledge of functional relationships (influences) between variables. This is the next major line of research that we hope to pursue.

## Appendix A: Derivations of key results

### A.1 Derivation of result (1)

An initial consideration is that the observations could have physical units (e.g. mass, length and time). MOEs will have to be stated consistently in some standard set of units. This is only possible if they are dimensionless. Otherwise, one would eventually meet cases where the required units could not be produced from the units available in the problem domain. Therefore, all MOEs must be dimensionless. To take matters further, we need to introduce the minimal algebra[1] for our purposes: $\langle (\cap, \cup, \neg), (+, -, \times, \div), (V_o, V_s) \rangle$. It includes the set theory operators, the arithmetic operators and, finally, the only objects for this problem, $\left(V_o, V_s\right)$. It is convenient to remove the complementation operator $\neg$, but to expand the objects to include their complements. This produces the algebra: $\left\langle \left(\cap, \cup\right), \left(+, -, \times, \div\right), \left(V_o, \overline{V}_o, V_s, \overline{V}_s\right) \right\rangle$, which can be shown to be equivalent to the first through de Morgan's laws. All of the elements can be treated as symbols, which can be arranged in various sequences to provide candidate MOE formulae. The sequences are subject to certain selection rules. With the latter algebra, one cannot have unbroken sequences of objects or of operators. Also, any sequence must begin and end with an object. These rules allow one to define a complexity measure for a mathematical definition of a MOE as the number of operators appearing in it (or as the number of objects – generally greater by unity). Brackets to clarify the scope of operators might be included in the algebra, but these are not needed to find that the simplest definition of a MOE capable of satisfying the requirements is given by result (1). This is a matter of generating increasingly complex legal sequences and testing them against the requirements. To illustrate the process, some unsuccessful MOE candidates are explained as follows:

- $V_o$ - positive, but no account taken of $V_s$, and could be dimensional.
- $V_o / V_s$ - positive and dimensionless, but could take infinite value.
- $\overline{V}_o \cap V_s / \overline{V}_o$ - dimensionless and values in range [0,1], but value decreases as acceptability of observations increases.

---

[1] An algebra is a set of objects, e.g. sets, and operations on those objects, e.g. union and intersection.

### A.2 Derivation of result (2)

With reference to Figure 2, which shows a 2-D example, each volume (observation or acceptance) can be decomposed into shells containing observations with equal probability, i.e. iso-probability shells. (Shells are surfaces with a finite thickness.) Consider shells labelled $i$ and $j$, respectively for acceptability and observations. The generalisation of the numerator of (1) is:

$$Q = \sum_{i=1}^{I} \sum_{j=1}^{J} F\left(\rho_{si}, \rho_{oj}\right) v_{si} \cap v_{oj} \qquad \text{(A-1)}$$

where $F\left(\rho_{si}, \rho_{oj}\right)$ is some function of the shell probability densities (PDs), and $v_{si}$ and $v_{oj}$ are the shell finite volumes. For acceptability and observations, different iso-probability shells are disjoint. Therefore, the intersections $v_{si} \cap v_{oj}$ are unique (no other pairs of shells have the same intersection) and complete (they cover all the space). If we let the shells become infinitesimal and their number infinite, such that all the space is still covered, result (A-1) becomes the integral:

$$Q = \int F\left(\rho_s(\underline{x}), \rho_o(\underline{x})\right) d\underline{x} \qquad \text{(A-2)}$$

where $\underline{x}$ is a position vector in the space, $\rho_s(\underline{x})$ and $\rho_o(\underline{x})$ are probability density functions (PDFs), and the integration is over all space. Note that the PDFs in (A-2) must be treated symmetrically because we are comparing like with like.

The 2-D example just treated involves shell intersections, which are immediately finite volume elements. This is not so in higher dimensions where the intersections may extend considerably through the space (e.g. in 3-D, intersecting spherical surfaces may produce circles). However, these intersections can be further decomposed into finite volume elements and result (A-2) still stands.

The function $F$ must be the simplest that allows a general, usable definition of MOE in terms of the PDFs. It is sufficient to consider functions involving just the arithmetic operators $(+, -, \times, \div)$. To reflect the required symmetry of the PDF arguments, $F$ should be commutative. Therefore, we are restricted to the operators $(+, \times)$. The addition operator leads to:

$$Q_+ = \int \left( \rho_s(\underline{x}) + \rho_o(\underline{x}) \right) d\underline{x} = 2 \tag{A-3}$$

Since $Q_+$ is invariant, it cannot lead to a usable definition of MOE. However, the remaining multiplication operator does, via the PDF overlap integral:

$$Q_\times = \int \rho_s(\underline{x}) \, \rho_o(\underline{x}) \, d\underline{x} \tag{A-4}$$

This should not be surprising as, by probability theory, the intersection of independent events, which effectively is what we have in the numerator of result (1), is associated with a product of their probabilities. Quantity $Q_\times$ has a lower limit of zero and a variable upper limit. To become a MOE consistent with result (1), it must be normalised to have an upper limit of unity. The upper limit is reached when all the observation probability mass is located at $\hat{\underline{x}}$ corresponding to the global maximum value of $\rho_s(\underline{x})$, denoted by $\rho_{s\max}$. Then, we have:

$$Q_{\max} = \int \rho_s(\underline{x}) \, \delta(\underline{x} - \hat{\underline{x}}) \, d\underline{x} = \rho_s(\hat{\underline{x}}) = \rho_{s\max} \tag{A-5}$$

Therefore, the general MOE for probabilistic systems becomes:

$$M = \frac{Q_\times}{Q_{\max}} = \int \frac{\rho_s(\underline{x})}{\rho_{s\max}} \, \rho_o(\underline{x}) \, d\underline{x} \tag{A-6}$$

### A.3   Derivation of result (22)

Rather than give a formal, general derivation we illustrate how the combination functions are obtained by considering a sufficiently rich example. This involves three users whose acceptability PDFs are determined by the parameter values $\underline{S}_1$, $\underline{S}_2$ and $\underline{S}_3$. We wish to find a user function reflecting satisfaction of pairs of users. The corresponding acceptability PDF is:

$$\rho_s(\underline{x}) = \rho \left( \underline{x} \, | \, \left( \underline{S}_1 \cap \underline{S}_2 \right) \cup \left( \underline{S}_1 \cap \underline{S}_3 \right) \cup \left( \underline{S}_2 \cap \underline{S}_3 \right) \right) \tag{A-7}$$

By Bayes' theorem, it can be expressed as the proportionality:

$$\rho_s(\underline{x}) \propto \rho \left( \left( \underline{S}_1 \cap \underline{S}_2 \right) \cup \left( \underline{S}_1 \cap \underline{S}_3 \right) \cup \left( \underline{S}_2 \cap \underline{S}_3 \right) \, | \, \underline{x} \right) \rho(\underline{x}) \tag{A-8}$$

where the first PDF on the RHS is the likelihood function and the second is the prior. The prior is unknown, so a non-informative one is used. This amounts to replacing $\rho(\underline{x})$ by a constant, which can be absorbed into the proportionality constant.

Consider a standard probability result expressed for PDFs over multi-dimensional $\underline{y}$ space say:

$$\rho\left(\underline{y}_1 \cup \underline{y}_2\right) = \rho\left(\underline{y}_1\right) + \rho\left(\underline{y}_2\right) - \rho\left(\underline{y}_1 \cap \underline{y}_2\right) \tag{A-9}$$

When the $\underline{y}_1$ and $\underline{y}_2$ values differ in at least one coordinate, the intersection, $\underline{y}_1 \cap \underline{y}_2$, lies in a sub-space of $\underline{y}$. Consequently, the PDF of the intersection is negligible compared with the other PDFs, and result (A-9) reduces to:

$$\rho\left(\underline{y}_1 \cup \underline{y}_2\right) = \rho\left(\underline{y}_1\right) + \rho\left(\underline{y}_2\right) \tag{A-10}$$

A similar result can be used to develop the likelihood function in result (A-8). This will produce conditional PDFs of double intersections, e.g. $\rho\left(\underline{S}_1 \cap \underline{S}_2 \mid \underline{x}\right)$. We assume that the users and, therefore, their acceptance parameter values are independent. So, for example:

$$\rho\left(\underline{S}_1 \cap \underline{S}_2 \mid \underline{x}\right) = \rho\left(\underline{S}_1 \mid \underline{x}\right)\rho\left(\underline{S}_2 \mid \underline{x}\right) \tag{A-11}$$

Bayes' theorem is applied to the PDFs on the RHS of result (A-11) to yield:

$$\rho\left(\underline{S}_1 \cap \underline{S}_2 \mid \underline{x}\right) \propto \rho\left(\underline{x} \mid \underline{S}_1\right)\rho\left(\underline{x} \mid \underline{S}_2\right)\rho\left(\underline{S}_1\right)\rho\left(\underline{S}_2\right) \tag{A-12}$$

Using results (A-10) to (A-12) in (A-8) produces:

$$\rho_s(\underline{x}) \propto \rho\left(\underline{x} \mid \underline{S}_1\right)\rho\left(\underline{x} \mid \underline{S}_2\right)\rho\left(\underline{S}_1\right)\rho\left(\underline{S}_2\right) + \rho\left(\underline{x} \mid \underline{S}_1\right)\rho\left(\underline{x} \mid \underline{S}_3\right)\rho\left(\underline{S}_1\right)\rho\left(\underline{S}_3\right) \\ + \rho\left(\underline{x} \mid \underline{S}_2\right)\rho\left(\underline{x} \mid \underline{S}_3\right)\rho\left(\underline{S}_2\right)\rho\left(\underline{S}_3\right) \tag{A-13}$$

The user function is initially taken as the RHS of result (A-13) divided by its maximum value. Provided all the user functions have their maxima at the same value of $\underline{x} = \hat{\underline{x}}$, the result can be written in terms of just user functions as follows:

$$f_s(\underline{x}) = W_{12}\ f_{s1}(\underline{x})\ f_{s2}(\underline{x}) + W_{13}\ f_{s1}(\underline{x})\ f_{s3}(\underline{x}) + W_{23}\ f_{s2}(\underline{x})\ f_{s3}(\underline{x}) \tag{A-14}$$

where the weights are given by:

$$W_{ij} = \frac{\rho\left(\underline{S}_i\right)\rho\left(\underline{S}_j\right)\rho\left(\hat{\underline{x}} \mid \underline{S}_i\right)\rho\left(\hat{\underline{x}} \mid \underline{S}_j\right)}{\sum_{ij} \rho\left(\underline{S}_i\right)\rho\left(\underline{S}_j\right)\rho\left(\hat{\underline{x}} \mid \underline{S}_i\right)\rho\left(\hat{\underline{x}} \mid \underline{S}_j\right)} \qquad \text{for} \qquad ij = 12, 13, 23 \tag{A-15}$$

The weights are positive and sum to unity. In practice, they would not be calculated using result (A-15), because the various PDFs featured are not available. They would be chosen to reflect the perceived relative importance of pairs of users. Result (A-15) merely serves to show the existence of weights that can be chosen at will. If there are no preferences between pairs of users, the weights are identically equal to the inverse of the number of pairs, i.e. $W_{ij} = 1/3$ for all $ij$.

It will be noticed that result (A-14) has the form of a user function squared. This is easily seen by making all the user functions equal, i.e. $f_{si}(\underline{x}) = f_s^{'}(\underline{x})$ for all '$i$' producing: $f_s(\underline{x}) = \left(f_s^{'}(\underline{x})\right)^2$. To obtain a simple user function, the square root of the latter must be taken. This also ensures that when all the user functions are equal, any combination of them is equal to just one of them. The final result is:

$$f_s(\underline{x}) = \sqrt{W_{12}\ f_{s1}(\underline{x})\ f_{s2}(\underline{x}) + W_{13}\ f_{s1}(\underline{x})\ f_{s3}(\underline{x}) + W_{23}\ f_{s2}(\underline{x})\ f_{s3}(\underline{x})} \tag{A-16}$$

This example demonstrates all the basic principles required to obtain the more complicated result (23) and the general result (22). The binomial factor in result (22) is the number of possible (distinguishable) $k$-tuplets with $I$ users.

### A.4 Derivation of result (28)

The error of the observation from the system being evaluated is $\underline{\varepsilon} = \underline{x} - \underline{z}$, where $\underline{x}$ is the observation and $\underline{z}$ is the corresponding ground truth. The error of the observation from the more accurate system is $\underline{\lambda} = \underline{y} - \underline{z}$ where $\underline{y}$ is the observation. Eliminating the unknown $\underline{z}$ from these two relations produces $\underline{\varepsilon} = \underline{x} - \underline{y} + \underline{\lambda}$. In terms of other PDFs, the PDF of the observation error is given by:

$$\rho_o(\underline{\varepsilon}) = \iiint \delta\left(\underline{\varepsilon} - \underline{x} + \underline{y} - \underline{\lambda}\right) \rho_x(\underline{x})\ \rho(\underline{y}, \underline{\lambda})\, d\underline{x}\ d\underline{y}\ d\underline{\lambda} \tag{A-17}$$

The PDFs in the integral allow for independence between $\underline{x}$ and $\underline{y}$ observations, and dependence between $\underline{y}$ and $\underline{\lambda}$. For single sample observations $\underline{x}^{'}$ and $\underline{y}^{'}$, the PDFs become:

$$\rho_x(\underline{x}) = \delta\left(\underline{x} - \underline{x}^{'}\right) \qquad \text{and} \qquad \rho\left(\underline{y}, \underline{\lambda}\right) = \delta\left(\underline{y} - \underline{y}^{'}\right) \rho_\lambda\left(\underline{\lambda}\right) \tag{A-18}$$

Inserting these relations into result (A-17) and carrying out the integrations over $\underline{x}$ and $\underline{y}$, and then $\underline{\lambda}$ yields:

$$\rho_o(\underline{\varepsilon}) = \rho_\lambda\left(\underline{\varepsilon} - \underline{x}^{'} + \underline{y}^{'}\right) \tag{A-19}$$

Inserting this into the usual MOE formula yields, finally:

$$M = \int f_s(\underline{\varepsilon})\ \rho_o(\underline{\varepsilon})\, d\underline{\varepsilon} = \int f_s(\underline{\varepsilon})\ \rho_\lambda\left(\underline{\varepsilon} - \underline{x}^{'} + \underline{y}^{'}\right) d\underline{\varepsilon} \tag{A-20}$$